\documentstyle[10pt]{amsart}

\title[Theorems of Barth-Lefschetz type]{Theorems of Barth-Lefschetz type
on K\"ahler manifolds of non-negative
bisectional curvature}
\author{Meeyoung Kim\and Jon Wolfson}
\thanks{The second author was partially supported by NSF grant DMS-9802487.}
\date{\today}

\newtheorem{thm}{Theorem}[section]
\newtheorem{lem}[thm]{Lemma}
\newtheorem{cor}[thm]{Corollary}
\newtheorem{prop}[thm]{Proposition}

\theoremstyle{definition}

\newtheorem{defn}{Definition}[section]

\numberwithin{equation}{section}

\renewcommand{\a}{\alpha}
\renewcommand{\b}{\beta}
\renewcommand{\d}{\delta}

\newcommand{\g}{\gamma}

\renewcommand{\l}{\lambda}

\newcommand{\n}{\nabla}

\newcommand{\Sig}{\Sigma}

\renewcommand{\O}{\Omega}
\renewcommand{\o}{\omega}

\def\Pb{\ifmmode{\Bbb P}\else{$\Bbb P$}\fi}
\def\Z{\ifmmode{\Bbb Z}\else{$\Bbb Z$}\fi}
\def\C{\ifmmode{\Bbb C}\else{$\Bbb C$}\fi}
\def\R{\ifmmode{\Bbb R}\else{$\Bbb R$}\fi}

\def\Ca{\cal C}
\newcommand\pn[1]{{\Bbb P}^{#1}}

\begin{document}

\maketitle

\setcounter{section}{-1}

\section{\bf Introduction}

The general philosophy  of Barth-Lefschetz type theorems is that
a subvariety $Y$ of ``small" codimension in a fixed variety $X$ must be
subject to stringent
topological restrictions.
For example, in the case of a complete intersection $Y$ in the
$n$-dimensional complex projective space $\pn{n}$,
the Lefschetz hyperplane theorem gives  cohomological restrictions on $Y$.

In 1970 W. Barth [B] discovered that even when $Y$ is not a complete
intersection, one can compare the cohomology of $\pn{n}$ to that
of $Y$ and prove theorems similar to  the Lefschetz theorem.
Since Barth's foundational work, there have been many studies in this
direction,
for example, those by Sommese [S] and Fulton-Lazarsfeld [F-L].

In [S-W], R. Schoen and J. Wolfson showed that a variant of an argument of
Frankel
together with Morse theory on a space of paths leads to a proof of
homotopic Barth-Lefschetz type theorems for complex submanifolds of
compact K\"ahler manifolds of non-negative holomorphic bisectional curvature.
To use the technique of [S-W] in manifolds other than $\pn{n}$  it is
necessary to compute a numerical invariant,
called the {\it complex positivity}, which
measures the positivity of the holomorphic bisectional curvature.
In this note, we compute the complex positivity  of all the compact
hermitian symmetric spaces. In fact this number plays an important role
in much work centered on hermitian symmetric spaces, for example the
metric rigidity theorems of Mok [M1].
Combining the computations of complex positivity with the results of [S-W]
we conclude:

\begin{thm}
\label{thm:intro}
Let $V$ be a K\"ahler manifold. Suppose that $M, N\subset V$ are 
compact complex
      submanifolds of complex dimensions $m,n$, respectively.
      Let
      $$i_*: \pi_j (N,N\cap M)\rightarrow \pi_j(V, M),$$
      be the homomorphism induced by the inclusion.
\begin{enumerate}
\item[(i)] If $V=U(p+q)/S(U(p) \times U(q))$, the complex grassmannian,
then $i_*$ is an isomorphism for $j\le n+m-2pq+p+q-1,$ and is surjective
for $j= n+m-2pq+p+q.$
\item[(ii)] If $V=SO(2 + p)/SO(2) \times SO(p)$, the complex quadric,
then $i_*$ is an isomorphism for $j\le n+m-p-1,$ and is surjective for
$j= n+m-p.$
\item[(iii)] If $V=Sp(r)/U(r),$
then  $i_*$ is an isomorphism for $j\le n+m-r^2,$ and is surjective
for $j= n+m-r^2+1.$
\item[(iv)] If $V=SO(2r)/U(r),$
then  $i_*$ is an isomorphism for $j\le n+m-r^2 +3r-3,$ and is surjective for
$j= n+m-r^2 +3r-2.$
\item[(v)] If $V=E_{6}/Spin(10)\times T^{1},$ then
$i_*$ is an isomorphism for $j\le n+m-21,$ and is surjective for
$j= n+m-20.$
\item[(vi)] If $V=E_{7}/E_{6}\times T^{1},$ then
$i_*$ is an isomorphism for $j\le n+m-37,$ and is surjective for
$j= n+m-36.$
\end{enumerate}
\end{thm}

Case (i) is stated in [S-W] and case (ii) is stated incorrectly
there, both without
proof. We remark that this result also follows from
computations of Goldstein [G] and Sommese's approach to Barth-Lefschetz
type theorems. On the other hand using the result of [S-W] the proof of
the theorem is reduced to what are essentially standard computations
in hermitian symmetric spaces. Further, the method of proof
of the theorem yields better results when some additional information
is available about the ambient curvature along either $M$ or $N$.
This should prove especially useful if either $M$ or $N$ is a hypersurface.
This is discussed in more detail
in Section 1. Finally the method described in [S-W] requires that
the ambient manifold be K\"ahler of non-negative biholomorphic sectional
curvature. Up to taking covers, by a uniformation theorem of Mok, every such
manifold is a product of flat space and hermitian symmetric spaces. 
Consequently
it was necessary to understand how well the method of [S-W] worked in
hermitian symmetric spaces. This is accomplished here.

\medskip
In June 2001 Meeyoung Kim died in a tragic drowning accident. 
She was 35 years old and had a promising mathematical future ahead of her.
Her mathematics was centered around the relation between topology and
algebraic geometry.
She had wanted the work reported on in this paper to be a beginning,
not an end. Her passing leaves all who knew
her saddened with a profound sense of loss.

\bigskip

\section{ \bf The index of a critical point}

Let $V$ be a  complete K\"ahler manifold of complex dimension $v$, with complex
structure $J$ and Levi-Civita connection $\n$.  Let $M$ and $N$ be
complex submanifolds of complex dimensions $m$ and $n$, respectively.
We denote, by ${\O}(V;M,N) = {\O}$, the space
of piecewise smooth paths $\g:[0,1] \to V$ constrained by the
requirements that $\g(0) \in M$ and $\g(1) \in N$.
Consider the energy of a path
$$ E(\g) = \int^1_0 |\dot{\g}|^2 dt $$
as a function on ${\O}$.  It is  shown in [S-W] that
$\g$ is a critical point of $E$ if:
\begin{enumerate}
\item[(i)] $\g$ is a smooth geodesic
\item[(ii)] $\g$ is normal to $M$ and $N$ at $\g(0)$ and $\g(1)$,
respectively.
\end{enumerate}
Let $W_1,W_2 \in T_{\g} {\O}$.  If $\g$ is a critical point of $E$ then
the second variation of $E$ along $\g$  is:
\begin{align}
\frac12 E_{**}(W_1,W_2) & = \langle \n_{W_1} W_2, \dot{\g}\rangle
\Big\vert^{1}_{0} + \int^1_0 \langle \n_{\dot{\g}} W_1, \n_{\dot{\g}}W_2
\rangle dt  \nonumber \\
& - \int^1_0 \langle R(\dot{\g}, W_1) \dot{\g}, W_2 \rangle dt.
\end{align}
where $R$ denotes the curvature tensor of $V$.

Suppose that $\g$ is a nontrivial critical point  and that $W(0)$ is a vector
in $T_{\g(0)} M$.  Parallel translate $W(0)$ along $\g$ to construct a vector
field $W = W(t)$ along $\g$.  Of course, $W(1)$ need not be tangent to $N$
at $\g(1)$ so $W$ is not necessarily an element of $T_{\g}{\O}$. However
formally we have:
\begin{equation}
\frac12 E_{**}(W,W) = \langle \n_{W} W, \dot{\g}\rangle
\Big\vert^{1}_{0} - \int^1_0 \langle  R(\dot{\g}, W) \dot{\g}, W \rangle
dt.
\label{equ:WW}
\end{equation}
$V$ is K\"ahler so $JW$ is also parallel along $\g$. $M$ is complex so
$JW(0) \in T_{\g(0)} M$.
Thus we also have:
\begin{equation}
\frac12 E_{**}(JW,JW) = \langle \n_{JW} JW, \dot{\g}\rangle
\Big\vert^{1}_{0} - \int^1_0 \langle  R(\dot{\g}, JW) \dot{\g}, JW \rangle
dt.
\label{equ:JJ}
\end{equation}
Adding (\ref{equ:WW}) and (\ref{equ:JJ}) and using $\n_{JW} JW = - \n_W W$
we have:
\begin{align}
\frac12 E_{**}(W,W) & + \frac12 E_{**} (JW,JW) =  \nonumber \\
& - \int^1_0(\langle R(\dot{\g},W)\dot{\g},W\rangle + \langle R(\dot{\g}, JW)
\dot{\g}, JW \rangle)dt.
\label{equ:sum}
\end{align}
Using the symmetries of the curvature tensor we have:
\begin{equation}
\langle R(\dot{\g},W)\dot{\g},W\rangle + \langle
R(\dot{\g},JW)\dot{\g},JW\rangle = \langle R(\dot{\g},J
\dot{\g})W,JW \rangle.
\label{equ:curv}
\end{equation}
This expression is the holomorphic bisectional curvature of the complex
lines $\dot{\g} \wedge J \dot{\g}$ and $W \wedge JW$.

Let $\{ W_1(0), \dots, W_m(0), JW_1 (0), \dots, JW_m(0)\}$ be an
orthonormal framing of $T_{\g(0)} M$.  For each $i = 1, \dots, m$, parallel
translate $W_i(0)$ along $\g$ to construct parallel vector fields
$\{ W_1, \dots, W_m, JW_1 , \dots, JW_m\}$ along $\g$.
Note that the vectors $W_i(1), JW_i(1)$ are perpendicular to both
$\dot{\g}(1)$ and $J\dot{\g}(1)$.  Thus the vector space
$$ S = \mbox{span}\{ W_1(1), \dots, W_m(1), JW_1(1), \dots, JW_m(1)\}$$
is a complex $m$-dimensional space lying in a complex $(v-1)$-dimensional
subspace of $T_{\g(1)}V$.  It follows that the subspace $S \cap  T_{\g(1)}
N$ has complex dimension at least equal to
$$ m+n - (v-1). $$
Moreover, the vector fields $\{W,JW\}$ with $W(1), JW(1) \in S \cap
T_{\g(1)}N$ are parallel and lie in $T_{\g} {\O}$.

\bigskip

\begin{thm}
\label{thm:positive}
Suppose that $V$ is a K\"ahler manifold of positive
holomorphic bisectional curvature, that $M$ and $N$ are complex submanifolds
and that $\g$ is a nontrivial critical point of the energy on
${\O}(V;M,N)$.  Then,
$$ \mbox{index} (\g) \ge  m+n-(v-1). $$
\end{thm}

\bigskip

\begin{pf}
There are at least $m+n-(v-1)$ pairs $\{W,JW\}$ that are parallel
along $\g$ and lie in $T_{\g}{\O}$. For each such pair, using the curvature
assumption (\ref{equ:sum}) and (\ref{equ:curv}) we have:
\begin{align*}
      & E_{**}(W,W)  +  E_{**} (JW,JW) =   \\
& - 2\int^1_0 \langle R(\dot{\g},J\dot{\g})W,JW \rangle  dt < 0.
\end{align*}
The result follows.
\end{pf}

\bigskip

Let $V$ be a K\"ahler manifold.
Fix $x \in V$ and let $X \wedge JX$ be a complex line in $T_x V$.  Let
${\Ca}(x, X
\wedge JX)$ be the cone:
$$ {\Ca}(x, X \wedge JX) = \{ Y \in T_x V: \langle R(X,JX) Y, JY \rangle >
0 \}. $$
Note that ${\Ca}$ is a complex cone; if $Y \in {\Ca}$ then $JY \in {\Ca}$.

\medskip

\begin{defn}
         \label{defn:positivity}
Denote the set of complex subspaces of ${\Ca}(x, X \wedge JX)$
by ${\cal L} ={\cal L} (x, X \wedge JX)$ and define
$$\ell (x,X \wedge JX) = \max_{L \in {\cal L}} \mbox{dim}_\C({ L}).$$
Then define,
\begin{enumerate}
\item[(i)] $\displaystyle{\ell(x) = \inf_{X \wedge JX} ~ \ell (x,X \wedge
JX)} $
\item[(ii)] $\displaystyle{\ell = \inf_{x \in V} ~ \ell(x)}.$
\end{enumerate}
We say that $\ell$ is the {\it complex positivity} of $V$.
\end{defn}
If $V$ is a hermitian symmetric space then clearly ${\Ca}(x, X \wedge JX)$
and hence $\ell (x,X \wedge JX)$ are independent of $x \in V$.

\medskip

\noindent{\bf Remark:} If $V$ is a compact hermitian symmetric space and $e$ is
the identity define the symmetric bilinear form $H_X(W, Z) = \langle
R(X,JX) W, JZ \rangle$,
where $X,W,Z \in T_eV$. Then for any $X \neq 0$, $H_X$ is positive
semi-definite.
Denote the null-space by ${\cal N}_X$. Then it follows easily that
$\ell (X \wedge JX)$
is the complementary dimension of  ${\cal N}_X$. For the four
classical families
of hermitian symmetric spaces it is easy to compute $\mbox{dim}_{\C}{\cal N}_X$
directly (see [M2]). For example, if $V = \mbox{Gr}(p,p + q;\C),$ the complex
Grassmann manifold and $X \in T_e V$ then $X \neq 0$ is a matrix with
$1 \leq \mbox{rank}(X) \leq \min(p,q)$. Then,
$$ \mbox{dim}_{\C}{\cal N}_X = (p - \mbox{rank}(X))(q - \mbox{rank}(X)).$$
Therefore, $\mbox{dim}_{\C}{\cal N}_X$ is maximal and $\ell(X \wedge JX)$
is minimal when $\mbox{rank}(X) = 1$.

\medskip

\begin{thm}
\label{thm:nonnegative}
\hskip .2cm
Suppose that $V$ is a complete K\"ahler manifold of non-negative holomorphic
bisectional curvature.
Let $M$ and $N$ be complex submanifolds of complex dimensions $m$ and $n$,
respectively,
and  $\g$ be a nontrivial
critical point of energy on ${\O}(V;M,N)$.  Then
$$ \mbox{index} (\g) \ge m + n - (v - 1) - (v - \ell). $$
\end{thm}

\begin{pf}
The argument in the proof of Theorem \ref{thm:positive} shows that if $W, JW
\in S \cap T_{\g(1)} N$ then
$$  E_{**}(W,W) + E_{**} (JW,JW) = - 2\int^1_0 \langle R(\dot{\g},
J\dot{\g}) W, JW \rangle \le 0. $$
To get strict inequality we want
$$ \langle R(\dot{\g},J(\dot{\g}))W,JW \rangle > 0$$
at  $\g(0)$.  This is insured by requiring  that:
$$W(0), JW(0) \in L \cap T_{\g(0)}M$$
where $L \in {\cal L}(\g(0),\dot{\g} \wedge J(\dot{\g})) $.
Let $L \in {\cal L}(\g(0),\dot{\g} \wedge J(\dot{\g}))$ be of maximal
dimension.
Then the complex dimension of $L \cap T_{\g(0)}M$ is at least
$\mbox{dim}_{\C} L + m -(v-1)$.
Parallel transport $L \cap T_{\g(0)}M$ along $\g$ to $\g(1)$ and denote the
resulting
subspace by $T$. Then $T \cap T_{\g(1)}N$ has complex dimension at least
$n + \mbox{dim}_{\C} L + m -(v-1) - (v-1)$.
The result follows.
\end{pf}

\begin{cor}
\label{cor:assumptionM}
Under the same hypotheses as the theorem assume that
$$ \ell(x, X \wedge JX) \geq \ell_0,$$
for every $x \in M$ and
every  complex line $X \wedge JX$ normal to $T_x M$ then
$$ \mbox{index} (\g) \ge m + n - (v - 1) - (v - \ell_0). $$
\end{cor}

\bigskip

\section{\bf{Applications}}

Let $V$ be a complete K\"ahler manifold of non-negative holomorphic bisectional
curvature, of complex dimension $v$ and with complex positivity $\ell$. Let
$M,N \subset V$ be
complex submanifolds of complex dimensions $m,n$, respectively and suppose that
$M$ is compact and $N$ is a closed subset of $V$.  The Morse theory of the
energy functional on
the path space $\O$ is described
in [S-W]. Combining this Morse theory with  Theorem \ref{thm:nonnegative}
it follows that:

\bigskip

\begin{thm}
\label{thm:homotopy}
Suppose that,
$$ \l_0 = n + m - v - (v - \ell) \geq 0. $$
Then relative homotopy groups $\pi_j(\O,N \cap M)$ are zero
for $0 \leq j \leq \l_0.$
\end{thm}

\bigskip

Theorem \ref{thm:homotopy} and the long exact homotopy sequence of the pair
$(\O, N \cap M)$
imply  that the homomorphism induced by the inclusion:
\begin{equation}
\imath_*:\pi_j(N \cap M) \to \pi_j(\O)
\label{equ:homotisom2}
\end{equation}
is an isomorphism  when $j < n+m - v - (v-\ell)$ and is a surjection when
$j = n+m - v - (v-\ell).$

Consider the fibration:
\begin{equation}
      \begin{array}{cccc}
{\O}(V;M,x)  \longrightarrow & {\O}(V;M,N)  \\
\\
& \downarrow e  \\
\\
& N
\end{array}
\end{equation}
where $e$ is the evaluation map $e: \g \mapsto \g(1)$ and  $x \in N$.
It is well-known that the homotopy groups of the fiber
${\O}(V;M,x)$ satisfy:
\begin{equation}
\pi_j ({\O}(V;M,x)) \simeq \pi_{j+1}(V,M),
\label{equ:homotisom1}
\end{equation}
for all $j$.  The long exact homotopy sequence of the fibration is:
\begin{align}
\cdots \longrightarrow \pi_{j+1}(N) & \longrightarrow \pi_j({\O}(V;M,x))
\longrightarrow \pi_j({\O}) \nonumber \\[.5cm]
&\overset{e_*}{\longrightarrow} \pi_j(N) \longrightarrow \pi_{j-1}
({\O}(V;M,x)) \longrightarrow  \cdots
\label{equ:fibrsequ}
\end{align}
Thus, using (\ref{equ:homotisom1}),
the long exact sequence (\ref{equ:fibrsequ}) becomes:
\begin{equation}
\cdots \rightarrow \pi_{j+1}(N)  \rightarrow
\pi_{j+1}(V,M) \rightarrow \pi_j(\O) \rightarrow \pi_j(N)
\rightarrow \pi_j (V,M) \rightarrow \cdots
\label{equ:homotseq}
\end{equation}
We have:

\bigskip

\begin{thm}
\label{thm:main}
Let $V$ be a complete K\"ahler manifold of non-negative holomorphic
bisectional curvature, of complex dimension $v$ and with complex positivity
$\ell$.   Let $M,N \subset V$
be complex submanifolds of complex dimensions $m,n$, respectively,  such that
$M$ is compact and  $N$ is a closed subset of $V$.  Then
the homomorphism induced by the inclusion
$$ \imath_*: \pi_j(N,N \cap M) \to \pi_j(V,M) $$
is an isomorphism for $j \le n+m-v-(v-\ell)$ and is a surjection for
$j = n+m-v-(v-\ell) + 1$ .
\end{thm}

\bigskip

\begin{pf}
For $\l_0 =  n+m-v-(v-\ell)$ consider the diagram:
\begin{equation*}
      \begin{array}{ccccc}
\pi_{{\l_0} +1}(N)\!\!\!\! &
\rightarrow \pi_{{\l_0}+1}(V,M)\!\!\!\!\!  &
\rightarrow \pi_{\l_0}(\O)\!\!\!\!\!
&\rightarrow \pi_{\l_0}(N)\!\!\!\! &
\rightarrow \pi_{\l_0}(V,M)\\
\\
\uparrow \simeq  & \uparrow & \uparrow \mbox{onto} & \uparrow \simeq &
\uparrow \\
\\
\pi_{{\l_0}+1}(N)\!\!\!\! & \rightarrow\!\! \pi_{{\l_0}+1}(N,N \cap M)
\!\!\!\!\! &
\rightarrow\!\! \pi_{\l_0}(N \cap M)\!\!\!\!\! & \rightarrow\!\! \pi_{\l_0}(N)
\!\!\!\! &
\rightarrow\!\! \pi_{\l_0}(N,N \cap M)
\end{array}
\end{equation*}
The vertical arrows are induced by inclusion.  The top row is the long
exact sequence (\ref{equ:homotseq}).  The bottom row is the long exact
sequence of the
pair $(N,N \cap M)$.  The result follows using Theorem \ref{thm:homotopy}
and the commutativity of the diagram.
\end{pf}

\bigskip

\begin{cor}
\label{cor:main}
      Under the same hypotheses as in  Theorem \ref{thm:main}, if
$$j \le 2m-v-(v-\ell)+1$$
then
$$ \pi_j(V,M) = 0. $$
\end{cor}

\begin{pf}
Apply Theorem \ref{thm:main} to the case $N = M$.
\end{pf}

\bigskip

\begin{cor}
Under the same hypothesis as in Theorem \ref{thm:main}, if
$$j \le \min(2m-v-(v-\ell)+1, ~ n+m -v-(v-\ell))$$
then
$$ \pi_j(N,N \cap M) = 0. $$
\end{cor}

\begin{pf}
      Follows from Corollary \ref{cor:main} and  Theorem \ref{thm:main}.
\end{pf}

\bigskip

\noindent{\bf Remark:} Under the same hypotheses as Theorem
\ref{thm:main} assume, in addition,  that
$$ \ell(x, X \wedge JX) \geq \ell_0,$$
for every $x \in M$ and
every  complex line $X \wedge JX$ normal to $T_x M$. Then by Corollary
\ref{cor:assumptionM}
the ranges of validity of the results of this section are improved by
replacing $\ell$ by $\ell_0$. For example suppose $V = \mbox{Gr}(p,p
+ q;\C)$ and $M$
is a complex hypersurface. If every normal $(1,0)$-vector $X$ along
$M$ has rank$(X) \geq r > 1$
then for {\it any} complex submanifold $N$ the index range in Theorem \ref{thm:main}
is increased by $r-1$.

\bigskip

\section{\bf The Compact Hermitian Symmetric Spaces}

\bigskip

We exploit the curvature computation in [Bo], to compute the complex
positivity of the  hermitian symmetric spaces. Recall the notation of [Bo].
Let $G$ be a compact
connected simple Lie group and $K$ be the  identity component of the
fixed point
set of an involutive  automorphism of $G$. We
assume that $K$ has infinite center.
Let ${\frak g}$ and ${\frak k}$ be the Lie algebras of
$G$ and $ K$ and ${\frak p}$ be the orthogonal
complement of ${\frak k}$ in ${\frak g}$ with respect to the
Killing form. Let ${\frak t}$ be a Cartan subalgebra of ${\frak k}$.
Denote the complexifications of these Lie algebras by
${\frak g}_{\C}, {\frak k}_{\C}$ and ${\frak t}_{\C}$.
Let $\Sig$ be the system of roots of ${\frak g}$ with respect to
${\frak t}$. Let ${\frak b}_{\a} $
be the one-dimensional eigenspace of ${\frak t}_{\C}$ in ${\frak g}_{\C}$
corresponding to $\a\in \Sig$.
We have,
\begin{eqnarray}
\label{equ:eigen}
[{\frak b}_{\a}, {\frak b}_{\b}] & = & \left\{ \begin{array}{ll}
{\frak b}_{\a+\b},
&\mbox{ if $  \a+\b$  is a root}  \; (\a+\b\ne 0),\\
       0,
&\mbox{ if $ \a+\b$  is not a root }\;  (\a+\b\ne 0),
\end{array}
\right.\\
h_{\a} &\in & [{\frak b}_{\a}, {\frak b}_{-\a}]. \nonumber
\end{eqnarray}
There are elements $e_{\a}\in {\frak b}_{\a}$
such that
$$[h, e_{\a}]=2 \pi i\, \a (h) \, e_{\a},\quad [e_{\a}, e_{-\a}]
=(2\pi i)^{-1}\, h_{\a}, \quad
(h\in {\frak t}) .$$
$\frak g$ is spanned by $\frak t$ and $e_{\a}+ e_{-\a},\, i(e_{\a}- e_{-\a}),$
for $\a\in \Sig$.
Define $N_{\a, \b}$ by:
$$[e_\a, e_\b]= N_{\a, \b}\, e_{\a+\b},\qquad \a,\b\in \Sigma,\, \a+\b\ne
0.$$ Then ([Bo]),
\begin{eqnarray}
\label{equ:NN}
N_{\a,\b}=-N_{-\a,-\b}.
\end{eqnarray}
The set of positive complementary roots of $G/K$,
i.e. the set of positive roots of ${\frak g}$ which are not
the roots of ${\frak k}$, will be denoted by $\Psi$. We put
$${\frak n}^{+}=\sum_{\a\in \Psi} {\frak b}_{\a},\quad
{\frak n}^{-}=\sum_{-\a\in \Psi}{\frak b}_{\a}.$$
We identify $\frak p$
with the tangent space of $G/K$ at $K$ (or at the identity $e$).
Then $T_*G/K \otimes \C$ is identified with ${\frak n}^{+} \oplus 
{\frak n}^{-}$.
Complex conjugation of ${\frak p} \otimes \C$ becomes
$e_\a \mapsto e_{-\a}$  and the assignment $e_\a \mapsto e_\a +
e_{-\a}$ defines
an isomorphism of complex vector spaces  between
${\frak n}^+$ and $\frak p$.
Here we endow $\frak p$ with the complex structure induced from that of
$G/K$.

For $\a\in \Sig$, let $\o_\a$  be the Maurer-Cartan form on $G_\C$
(the complex
Lie group with Lie algebra ${\frak g}_\C$)  defined by,
$$
\o_\a (e_\a ) = 1,\;\;\;
\quad \o_\a (e_\b ) = 0,\;\; \a \ne \b,\;\;\;
\o_\a (t) = 0,\;\; t \in {\frak t}_\C.
$$
This induces a K\"ahler metric on $G/K$ with orthonormal basis $\{
e_\a : \a\in \Psi\}$.

Note that $e_{\a}$ and $e_{-\a}$ are  $(1,0)$ and   $(0,1)$ vectors,
respectively. The curvature tensor with respect to the
basis $\{ e_\a : \a\in \Psi\}$ is
$$R(e_{-\g}, e_{\d}) e_{\a}= {\rm ad}\, [e_{\d}, e_{-\g}] (e_{\a})
=[[e_{\d}, e_{-\g}], e_{\a}]=\sum_{\b} R^{\b}_{\a\bar{\g}\d} e_{\b}.$$
Writing,
$$
e_{\a}={1\over \sqrt{2}} (X_\a - iJX_\a ),\quad e_{-\a}={1\over
\sqrt{2}} (X_\a + iJX_\a ),
$$
the real tangent vectors $\{ X_\a : \a\in \Psi\}$ to
$G/K$ at $K$ form a basis
for ${\frak p}$ and using the K\"ahler symmetries of curvature,
\begin{equation*}
          \langle R(X_\a, JX_\b) X_\g, JX_\d\rangle=
\,\,- ~ \langle R (e_{\a}, e_{-\b}) e_{\g},\, e_{-\d}\rangle.
\end{equation*}
For $\a,\b\in\Sigma$, we will denote the Killing form
$\kappa (h_{\a}, h_{\b})= {\rm Tr} ({\rm ad} \, h_{\a}\, {\rm ad}\, h_{\b})$
by $(\a,\b)$.

\medskip

\begin{lem} {\rm ([Bo] 1.4)}
         \label{lem:killing}
If $\a,\b\in\Psi,$  then  $\a+\b$
is not a root and therefore $ (\a,\b)\ge 0.$

\end{lem}

We use the following theorem to compute the holomorphic bisectional curvature:

\medskip

\begin{thm} {\rm ([Bo] 2.1,2.2)}
         \label{thm:Bo}
$$R^{\g}_{\b\bar{\d}\a}= \left\{ \begin{array}{ll}
                              0, &\mbox{ if $\b+\a\ne\g+\d$}\\
                          (\b, \d), &\mbox{ if $\b=\g,\, \a=\d$}\\
                        N_{\a, -\d} N_{\b, -\g}, &\mbox{ if $\b+\a=\g+\d,\,
\b\ne\g.$ }
		 \end{array}
\right.$$
\end{thm}

\bigskip

\noindent Hence,

\begin{lem}
         \label{lem:bisec}
\begin{eqnarray*}
       \lefteqn
         {\!\!\!\!\!\!\!\!\!\!
         \langle R (\sum_{\a} a_\a X_{\a},\, J\sum_{\a} a_\a X_{\a})
\sum_{\b} b_\b X_{\b},
\, J\sum_{\b} b_\b X_{\b}\rangle}\\
&=&  ~ \sum_{\a,\g} a_{\a}^2 b_{\g}^2\, (\a,\g)
+  \sum_{\a\ne \b, \g\ne \d, \; \a+\g=\b+\d}
a_{\a} a_{\b} b_{\g} b_{\d}\, N_{\a,-\b}N_{\g,-\d}.
\end{eqnarray*}
\end{lem}

\begin{pf}
For $\a=\b$ and $\g\ne \d,$ by Theorem \ref{thm:Bo},
$$
\langle R (X_{\a}, J X_{\a} ) X_{\g},\, J X_{\d}\rangle
=- ~ \langle R (e_{\a}, e_{-\a} ) e_{\g},\, e_{-\d}\rangle
= ~  R^{\d}_{\g\bar{\a}\a}=0.
$$
Similarly, for $\a\ne\b$ and $\g=\d,$
$$
\langle R (X_{\a}, J X_{\b} ) X_{\g},\, J X_{\g}\rangle
=0.
$$
Therefore,
\begin{eqnarray*}
&&  \hspace{-.5cm} \langle R (\sum_{\a} a_\a X_{\a}, J\sum_{\a} a_\a X_{\a})
\sum_{\b} b_\b X_{\b},
  J\sum_{\b} b_\b X_{\b}\rangle\\
&=& \hspace{-.3cm} \sum_{\a,\g}  \hspace{-.05cm} a_{\a}^2 b_{\g}^2
\langle R (X_{\a}, JX_{\a}) X_{\g}, JX_{\g}\rangle
+  \hspace{-.4cm} \sum_{\a\ne \b, \g\ne \d}  \hspace{-.2cm} a_{\a} 
a_{\b} b_{\g} b_{\d}
\langle R (X_{\a}, JX_{\b}) X_{\g}, JX_{\d}\rangle.
\end{eqnarray*}

Using the orthonormality of $\{ e_{\a} : \a\in\Psi\}$ and
Theorem \ref{thm:Bo},
we have
$$
\langle R (X_{\a}, J X_{\a} ) X_{\g},\, J X_{\g}\rangle
= ~ R^{\g}_{\g\bar{\a}\a}\langle e_{\g},e_{-\g}\rangle
= ~ (\a,\g).
$$
For $\a\ne\b, \g\ne\d$, by Theorem \ref{thm:Bo}, we have
$$
\langle R (X_{\a},\, JX_{\b}) X_{\g},\, JX_{\d}\rangle = 
R^{\d}_{\g\bar{\b}\a} =
\left\{
         \begin{array}{ll}
           N_{\a,-\b}N_{\g,-\d}
                         & \mbox{if $\a+\g=\b+\d$,  }\\
         0 & \mbox{otherwise.}
         \end{array}
\right.
$$
The lemma follows.
\end{pf}

\bigskip

Suppose that $e_{\a_j}, e_{-\a_j}$ for $\a_j \in \Psi, j=1,  \dots,
k,$ lie in an abelian
subspace of ${\frak p}$. Then
\begin{equation}
\label{eqn:Nvanishes}
N_{\a_i, -\a_j} = 0 ,\;\; i \ne j ,\;\; i,j=1,\dots, k.
\end{equation}
Let $X_{\a_i}$ be the vector associated to the root $\a_i$. From
Lemma \ref{lem:bisec} it follows that:
$$
      \langle R ( X_{\a_i},\, J X_{\a_j}) Y,\, JY\rangle = 0, \; i \ne j,
$$
for any vector $Y$.
This implies:

\begin{lem}
Suppose that $X = \sum_i^k a_i X_{\a_i}$, where $e_{\a_i}, e_{-\a_i}$, $i=1,
\dots, k,$ lie in an abelian subspace of ${\frak p}$.
Then for any $i=1, \dots, k,$
$$
\ell(X_{\a_i} \wedge JX_{\a_i}) \leq \ell(X \wedge JX).
$$
\end{lem}

\begin{pf}
     From (\ref{eqn:Nvanishes}) and Lemma \ref{lem:bisec}:
$$
      \langle R ( \sum_{i=1}^k a_i X_{\a_i},\, J \sum_{i=1}^k a_i X_{\a_i})
\sum_\b b_\b X_\b,\, J\sum_\b b_\b X_\b \rangle
= \sum_{i=1}^k \sum_{\b } a_i^2 b_\b^2 (\a_i, \b).
$$
Since $(\g, \d) \ge 0$ for $\g,\d \in \Psi$, the result follows.
\end{pf}

Under the action of the linear isotropy group any tangent vector of a
hermitian symmetric space
can be written as the sum of vectors all of which lie in a maximal
abelian subspace ${\frak a} \subset {\frak p}$. Therefore the infimum
$$\inf_{X\wedge JX} \ell( X\wedge JX)$$
is achieved when $X$ is spanned
by one base element, i.e. $X=a_{\a} X_{\a}$.
Fix $\a\in\Psi$ and set $\Psi_{\a}' :=\{\g\in\Psi :  (\a,\g)\ne 0 \}.$

\bigskip
\begin{prop}
         \label{prop:conealpha}
      For  fixed $X_{\a}$:

      \begin{enumerate}

          \medskip
          \item[(i)]
\,\, $\langle R ( X_{\a},\, J X_{\a})
\sum_{\b}\, b_\b X_{\b},
\, J\sum_{\b}\, b_\b X_{\b}\rangle
= \sum_{\b}\,  a_{\a}^2 b_{\b}^2 \, (\a,\b).$

\medskip
\item[(ii)]
The span of $\{ X_{\g}, J X_{\g} : \g\in\Psi_{\a}'\}$ is a maximal
subspace of the cone
${\cal C} ( X_{\a}\wedge JX_{\a})$.
\end{enumerate}
\end{prop}

\begin{pf}
(i) immediately follows from Lemma \ref{lem:bisec}.  (ii)
then follows from (i).
\end{pf}

\medskip
\begin{cor}
         \label{cor:ell}
For a compact hermitian symmetric space
the complex positivity  is  equal to
$$\ell = | \Psi_{\a}'|,$$
for any $\a\in \Psi$.
\end{cor}

\begin{pf}
By Proposition \ref{prop:conealpha}, it follows that $\ell
=\inf_{\a\in \Psi} | \Psi_{\a}'|$.
However, the elements of $\Psi$ are all highest weight vectors and
hence they are equivalent
under the action of the isotropy group. Therefore $| \Psi_{\a}'|$ is
independent of $\a \in \Psi$.
\end{pf}

\bigskip

\noindent By Corollary \ref{cor:ell} the complex positivities of the
compact hermitian symmetric spaces
are:

\medskip

\noindent{\bf $\mbox{Gr}(p,p + q;\C)$}: \hspace{.81cm} $\ell = p + q -1.$

\noindent{\bf $\mbox{Gr}(2,p + 2;\R)$}: \hspace{.87cm} $\ell = p -1.$

\noindent{\bf $Sp(r)/U(r)$}: \hspace{1.45cm} $\ell = r.$

\noindent{\bf $SO(2r)\!/\!U(r)$}: \hspace{1.35cm} $\ell = (r-1) + (r-2).$

\noindent{\bf $E_{6}/(Spin(10)\times T^{1})$}: \; $\ell = 11$

\noindent{\bf  $E_{7}/(E_{6}\times T^{1})$}: \hspace{1.15cm} $\ell= 17$.

\medskip

\noindent Combining this computation with Theorem \ref{thm:main} 
gives Theorem \ref{thm:intro}.

\end{document}